%
\catcode`\@=11
%
\font\elevenocm=bbold10 at 10pt
\font\tifont=cmbx10 at 12pt
\font\sevenocm=bbold10 at 8pt
\font\fiveocm=bbold10 at 7pt

\newfam\ocmfam
\textfont\ocmfam=\elevenocm  \scriptfont\ocmfam=\sevenocm
  \scriptscriptfont\ocmfam=\fiveocm
\def\ocm{\ifmmode\let\next\ocm@\else
 \def\next{\errmessage{Use \string\ocm\space only in math mode}}\fi\next}
\def\ocm@#1{{\ocm@@{#1}}}
\def\ocm@@#1{\fam\ocmfam#1}
\catcode`\@=12

\def\frac#1#2{{#1\over#2}} 
\def\PROB{{\ocm P}}

\def\IND#1{{\ocm 1}_{{\left[ #1 \right]}}}

\def\EXP{{\ocm E}}
\def\FR#1#2{\raise 3pt\hbox{$\scriptstyle #1$}\kern-2.1pt/\kern-2pt\lower1.8pt\hbox{$\scriptstyle #2$}}

\def\Var{{\ocm V}}

\def\isdef{\buildrel {\rm def} \over =}

\def\square{{\vcenter{\vbox{\hrule height .4pt
  \hbox{\vrule width .4pt height 7pt \kern7pt
  \vrule width .4pt}
  \hrule height .4pt}}}}

\def\tendsinlaw{\buildrel {\cal L} \over \to}

\def\inlaw{\buildrel {\cal L} \over =}

%
%
%
%
%
%
%
%

\def\mathbox#1#2#3{{\vcenter{\vbox{\hrule height #3%
  \hbox{\vrule width #3 height #2 \kern#1%
  \vrule width #3}%
  \hrule height #3}}}}%
\def\square{{\vcenter{\vbox{\hrule height .4pt%
  \hbox{\vrule width .4pt height 7pt \kern7pt%
  \vrule width .4pt}%
  \hrule height .4pt}}}}%
\def\emptybox#1#2#3{{\vbox{\hrule height #3%
  \hbox{\vrule width #3 height #2 \kern#1%
  \vrule width #3}%
  \hrule height #3}}}%
\def\figure#1\par{%
   \medskip%
   \centerline{\vbox{\hrule width 3.5in}}\smallskip%
   \centerline{Put Figure {#1} about here.} \smallskip%
   \centerline{\vbox{\hrule width 3.5in}}%
   \medskip}%

\def\temp{1.32}
\let\tempp=\relax
\expandafter\ifx\csname psboxversion\endcsname\relax
  \message{PSBOX(\temp) loading}
\else
    \ifdim\temp cm>\psboxversion cm
      \message{PSBOX(\temp) loading}
    \else
      \message{PSBOX(\psboxversion) is already loaded: I won't load
        PSBOX(\temp)!}
      \let\temp=\psboxversion
      \let\tempp=\endinput
    \fi
\fi
\tempp
\let\psboxversion=\temp
\catcode`\@=11
%
%
\def\execute#1{#1}
\def\psm@keother#1{\catcode`#112\relax}
\def\executeinspecs#1{%
\execute{\begingroup\let\do\psm@keother\dospecials\catcode`\^^M=9#1\endgroup}}
%
%
\def\psfortextures{
\def\PSspeci@l##1##2{%
\special{illustration ##1\space scaled ##2}%
}}
\def\psfordvitops{
\def\PSspeci@l##1##2{%
\special{dvitops: import ##1\space \the\drawingwd \the\drawinght}%
}}
\def\psfordvips{
\def\PSspeci@l##1##2{%
\d@my=0.1bp \d@mx=\drawingwd \divide\d@mx by\d@my%
\includegraphics{##1\space}%
}}
\def\psforoztex{
\def\PSspeci@l##1##2{%
\special{##1 \space
      ##2 1000 div dup scale
      \putsp@ce{\number-\psllx} \putsp@ce{\number-\pslly} translate
}%
}}
\def\putsp@ce#1{#1 }
\def\psfordvitps{
\def\psdimt@n@sp##1{\d@mx=##1\relax\edef\psn@sp{\number\d@mx}}
\def\PSspeci@l##1##2{%
\special{dvitps: Include0 "psfig.psr"}
\psdimt@n@sp{\drawingwd}
\special{dvitps: Literal "\psn@sp\space"}
\psdimt@n@sp{\drawinght}
\special{dvitps: Literal "\psn@sp\space"}
\psdimt@n@sp{\psllx bp}
\special{dvitps: Literal "\psn@sp\space"}
\psdimt@n@sp{\pslly bp}
\special{dvitps: Literal "\psn@sp\space"}
\psdimt@n@sp{\psurx bp}
\special{dvitps: Literal "\psn@sp\space"}
\psdimt@n@sp{\psury bp}
\special{dvitps: Literal "\psn@sp\space startTexFig\space"}
\special{dvitps: Include1 "##1"}
\special{dvitps: Literal "endTexFig\space"}
}}
\def\psfordvialw{
\def\PSspeci@l##1##2{
\special{language "PostScript",
position = "bottom left",
literal "  \psllx\space \pslly\space translate
  ##2 1000 div dup scale
  -\psllx\space -\pslly\space translate",
include "##1"}
}}
\def\psonlyboxes{
\def\PSspeci@l##1##2{%
\at(0cm;0cm){\boxit{\vbox to\drawinght
  {\vss
  \hbox to\drawingwd{\at(0cm;0cm){\hbox{(##1)}}\hss}
  }}}
}%
}
\def\psloc@lerr#1{%
\let\savedPSspeci@l=\PSspeci@l%
\def\PSspeci@l##1##2{%
\at(0cm;0cm){\boxit{\vbox to\drawinght
  {\vss
  \hbox to\drawingwd{\at(0cm;0cm){\hbox{(##1) #1}}\hss}
  }}}
\let\PSspeci@l=\savedPSspeci@l
}%
}
%
%
\newread\pst@mpin
\newdimen\drawinght\newdimen\drawingwd
\newdimen\psxoffset\newdimen\psyoffset
\newbox\drawingBox
\newif\ifNotB@undingBox
\newhelp\PShelp{Proceed: you'll have a 5cm square blank box instead of
your graphics (Jean Orloff).}
\def\@mpty{}
\def\s@tsize#1 #2 #3 #4\@ndsize{
  \def\psllx{#1}\def\pslly{#2}%
  \def\psurx{#3}\def\psury{#4}
  \ifx\psurx\@mpty\NotB@undingBoxtrue
  \else
    \drawinght=#4bp\advance\drawinght by-#2bp
    \drawingwd=#3bp\advance\drawingwd by-#1bp
  \fi
  }
\def\sc@nline#1:#2\@ndline{\edef\p@rameter{#1}\edef\v@lue{#2}}
\def\g@bblefirstblank#1#2:{\ifx#1 \else#1\fi#2}
\def\psm@keother#1{\catcode`#112\relax}
\def\execute#1{#1}
{\catcode`\%=12
\xdef\B@undingBox{
}   
\def\ReadPSize#1{
 \edef\PSfilename{#1}
 \openin\pst@mpin=#1\relax
 \ifeof\pst@mpin \errhelp=\PShelp
   \errmessage{I haven't found your postscript file (\PSfilename)}
   \psloc@lerr{was not found}
   \s@tsize 0 0 142 142\@ndsize
   \closein\pst@mpin
 \else
   \immediate\write\psbj@inaux{#1,}
   \loop
     \executeinspecs{\catcode`\ =10\global\read\pst@mpin to\n@xtline}
     \ifeof\pst@mpin
       \errhelp=\PShelp
       \errmessage{(\PSfilename) is not an Encapsulated PostScript File:
           I could not find any \B@undingBox: line.}
       \edef\v@lue{0 0 142 142:}
       \psloc@lerr{is not an EPSFile}
       \NotB@undingBoxfalse
     \else
       \expandafter\sc@nline\n@xtline:\@ndline
       \ifx\p@rameter\B@undingBox\NotB@undingBoxfalse
         \edef\t@mp{%
           \expandafter\g@bblefirstblank\v@lue\space\space\space}
         \expandafter\s@tsize\t@mp\@ndsize
       \else\NotB@undingBoxtrue
       \fi
     \fi
   \ifNotB@undingBox\repeat
   \closein\pst@mpin
 \fi
\message{#1}
}
%
%
\newcount\xscale \newcount\yscale \newdimen\pscm\pscm=1cm
\newdimen\d@mx \newdimen\d@my
\newdimen\pswdincr \newdimen\pshtincr
\let\ps@nnotation=\relax
\def\psboxto(#1;#2)#3{\vbox{
   \ReadPSize{#3}
   \advance\drawingwd by\pswdincr
   \advance\drawinght by\pshtincr
   \divide\drawingwd by 1000
   \divide\drawinght by 1000
   \d@mx=#1
   \ifdim\d@mx=0pt\xscale=1000
         \else \xscale=\d@mx \divide \xscale by \drawingwd\fi
   \d@my=#2
   \ifdim\d@my=0pt\yscale=1000
         \else \yscale=\d@my \divide \yscale by \drawinght\fi
   \ifnum\yscale=1000
         \else\ifnum\xscale=1000\xscale=\yscale
                    \else\ifnum\yscale<\xscale\xscale=\yscale\fi
              \fi
   \fi
   \divide \psxoffset by 1000\multiply\psxoffset by \xscale
   \divide \psyoffset by 1000\multiply\psyoffset by \xscale
   \global\divide\pscm by 1000
   \global\multiply\pscm by\xscale
   \multiply\drawingwd by\xscale \multiply\drawinght by\xscale
   \ifdim\d@mx=0pt\d@mx=\drawingwd\fi
   \ifdim\d@my=0pt\d@my=\drawinght\fi
   \message{scaled \the\xscale}
 \hbox to\d@mx{\hss\vbox to\d@my{\vss
   \global\setbox\drawingBox=\hbox to 0pt{\kern\psxoffset\vbox to 0pt{
      \kern-\psyoffset
      \PSspeci@l{\PSfilename}{\the\xscale}
      \vss}\hss\ps@nnotation}
   \global\ht\drawingBox=\the\drawinght
   \global\wd\drawingBox=\the\drawingwd
   \baselineskip=0pt
   \copy\drawingBox
 \vss}\hss}
  \global\psxoffset=0pt
  \global\psyoffset=0pt
  \global\pswdincr=0pt
  \global\pshtincr=0pt 
  \global\pscm=1cm
  \global\drawingwd=\drawingwd
  \global\drawinght=\drawinght
}}
%
%
\def\psboxscaled#1#2{\vbox{
  \ReadPSize{#2}
  \xscale=#1
  \message{scaled \the\xscale}
  \advance\drawingwd by\pswdincr\advance\drawinght by\pshtincr
  \divide\drawingwd by 1000\multiply\drawingwd by\xscale
  \divide\drawinght by 1000\multiply\drawinght by\xscale
  \divide \psxoffset by 1000\multiply\psxoffset by \xscale
  \divide \psyoffset by 1000\multiply\psyoffset by \xscale
  \global\divide\pscm by 1000
  \global\multiply\pscm by\xscale
  \global\setbox\drawingBox=\hbox to 0pt{\kern\psxoffset\vbox to 0pt{
     \kern-\psyoffset
     \PSspeci@l{\PSfilename}{\the\xscale}
     \vss}\hss\ps@nnotation}
  \global\ht\drawingBox=\the\drawinght
  \global\wd\drawingBox=\the\drawingwd
  \baselineskip=0pt
  \copy\drawingBox
  \global\psxoffset=0pt
  \global\psyoffset=0pt
  \global\pswdincr=0pt
  \global\pshtincr=0pt 
  \global\pscm=1cm
  \global\drawingwd=\drawingwd
  \global\drawinght=\drawinght
}}
%
\def\psbox#1{\psboxscaled{1000}{#1}}
%
%
%
\newif\ifn@teof\n@teoftrue
\newif\ifc@ntrolline
\newif\ifmatch
\newread\j@insplitin
\newwrite\j@insplitout
\newwrite\psbj@inaux
\immediate\openout\psbj@inaux=psbjoin.aux
\immediate\write\psbj@inaux{\string\joinfiles}
\immediate\write\psbj@inaux{\jobname,}
%
%
\immediate\let\oldinput=\input
\def\input#1 {
 \immediate\write\psbj@inaux{#1,}
 \oldinput #1 }
\def\empty{}
\def\setmatchif#1\contains#2{
  \def\match##1#2##2\endmatch{
    \def\tmp{##2}
    \ifx\empty\tmp
      \matchfalse
    \else
      \matchtrue
    \fi}
  \match#1#2\endmatch}
\def\warnopenout#1#2{
 \setmatchif{TrashMe,psbjoin.aux,psbjoint.tex}\contains{#2}
 \ifmatch
 \else
   \immediate\openin\pst@mpin=#2
   \ifeof\pst@mpin
     \else
     \errhelp{If the content of this file is so precious to you, abort (ie
press x or e) and rename it before retrying.}
     \errmessage{I'm just about to replace your file named #2}
   \fi
   \immediate\closein\pst@mpin
 \fi
 \message{#2}
 \immediate\openout#1=#2}
{
\catcode`\%=12
\gdef\splitfile#1 {
 \immediate\openin\j@insplitin=#1
 \message{Splitting file #1 into:}
 \warnopenout\j@insplitout{TrashMe}
 \loop
   \ifeof\j@insplitin
     \immediate\closein\j@insplitin\n@teoffalse
   \else
     \n@teoftrue
     \executeinspecs{\global\read\j@insplitin to\spl@tinline\expandafter
       \ch@ckbeginnewfile\spl@tinline
     \ifc@ntrolline
     \else
       \toks0=\expandafter{\spl@tinline}
       \immediate\write\j@insplitout{\the\toks0}
     \fi
   \fi
 \ifn@teof\repeat
 \immediate\closeout\j@insplitout}
\gdef\ch@ckbeginnewfile#1
 \def\t@mp{#1}
 \ifx\empty\t@mp
   \def\t@mp{#3}
   \ifx\empty\t@mp
     \global\c@ntrollinefalse
   \else
     \immediate\closeout\j@insplitout
     \warnopenout\j@insplitout{#2}
     \global\c@ntrollinetrue
   \fi
 \else
   \global\c@ntrollinefalse
 \fi}
\gdef\joinfiles#1\into#2 {
 \message{Joining following files into}
 \warnopenout\j@insplitout{#2}
 \message{:}
 {
 \edef\w@##1{\immediate\write\j@insplitout{##1}}
 \w@{
 \w@{
 \w@{
 \w@{
 \w@{
 \w@{
 \w@{
 \w@{
 \w@{\string\input\space psbox.tex}
 \w@{\string\splitfile{\string\jobname}}
 }
 \tre@tfilelist#1, \endtre@t
 \immediate\closeout\j@insplitout}
\gdef\tre@tfilelist#1, #2\endtre@t{
 \def\t@mp{#1}
 \ifx\empty\t@mp
   \else
   \llj@in{#1}
   \tre@tfilelist#2, \endtre@t
 \fi}
\gdef\llj@in#1{
 \immediate\openin\j@insplitin=#1
 \ifeof\j@insplitin
   \errmessage{I couldn't find file #1.}
   \else
   \message{#1}
   \toks0={
   \immediate\write\j@insplitout{\the\toks0}
   \executeinspecs{\global\read\j@insplitin to\oldj@ininline}
   \loop
     \ifeof\j@insplitin\immediate\closein\j@insplitin\n@teoffalse
       \else\n@teoftrue
       \executeinspecs{\global\read\j@insplitin to\j@ininline}
       \toks0=\expandafter{\oldj@ininline}
       \let\oldj@ininline=\j@ininline
       \immediate\write\j@insplitout{\the\toks0}
     \fi
   \ifn@teof
   \repeat
   \immediate\closein\j@insplitin
 \fi}
}
\def\autojoin{
 \immediate\write\psbj@inaux{\string\into\space psbjoint.tex}
 \immediate\closeout\psbj@inaux
 \input psbjoin.aux
}
%
%
%
%
\def\centinsert#1{\midinsert\line{\hss#1\hss}\endinsert}
\def\psannotate#1#2{\vbox{
  \def\ps@nnotation{#2\global\let\ps@nnotation=\relax}#1}}
\def\pscaption#1#2{\vbox{
   \setbox\drawingBox=#1
   \copy\drawingBox
   \vskip\baselineskip
   \vbox{\hsize=\wd\drawingBox\setbox0=\hbox{#2}
     \ifdim\wd0>\hsize
       \noindent\unhbox0\tolerance=5000
    \else\centerline{\box0}
    \fi
}}}
%
\def\at(#1;#2)#3{\setbox0=\hbox{#3}\ht0=0pt\dp0=0pt
  \rlap{\kern#1\vbox to0pt{\kern-#2\box0\vss}}}
%
\newdimen\gridht \newdimen\gridwd
\def\gridfill(#1;#2){
  \setbox0=\hbox to 1\pscm
  {\vrule height1\pscm width.4pt\leaders\hrule\hfill}
  \gridht=#1
  \divide\gridht by \ht0
  \multiply\gridht by \ht0
  \gridwd=#2
  \divide\gridwd by \wd0
  \multiply\gridwd by \wd0
  \advance \gridwd by \wd0
  \vbox to \gridht{\leaders\hbox to\gridwd{\leaders\box0\hfill}\vfill}}
%
\def\fillinggrid{\at(0cm;0cm){\vbox{
  \gridfill(\drawinght;\drawingwd)}}}
%
%
\def\textleftof#1:{
  \setbox1=#1
  \setbox0=\vbox\bgroup
    \advance\hsize by -\wd1 \advance\hsize by -2em}
\def\textrightof#1:{
  \setbox0=#1
  \setbox1=\vbox\bgroup
    \advance\hsize by -\wd0 \advance\hsize by -2em}
\def\endtext{
  \egroup
  \hbox to \hsize{\valign{\vfil##\vfil\cr%
\box0\cr%
\noalign{\hss}\box1\cr}}}
%
\def\frameit#1#2#3{\hbox{\vrule width#1\vbox{
  \hrule height#1\vskip#2\hbox{\hskip#2\vbox{#3}\hskip#2}%
        \vskip#2\hrule height#1}\vrule width#1}}
\def\boxit#1{\frameit{0.4pt}{0pt}{#1}}

\newcount\figurenumber
\figurenumber=0

\outer\def\psfigure#1#2%
{\par
\medskip
\vbox{\baselineskip=10.8pt%
\global\advance\figurenumber by 1
$$\pscaption{\psbox{#1}}
{\rm%
{\bf Figure \the\figurenumber.\enspace}
#2}$$
\par\medskip
}
}

\catcode`\@=12 
%
 \psfordvips   

\catcode`@=11 
\newskip\ttglue
\def\twelvepoint{%
\font\sc=cmcsc10
\font\twelverm=cmr10
\font\twelvei=cmmi10
\font\twelvebi=cmmib10

\font\twelveex=cmex10
\font\twelvesy=cmsy10
\font\twelvebf=cmbx10
\font\twelvett=cmtt10
\font\twelvesl=cmsl10
\font\twelvesc=cmcsc10
\font\twelveit=cmti10
\font\twelveu=cmu10
\font\ninerm=cmr10 
\font\eightrm=cmr9 at 8pt
\font\eighti=cmmi9 at 8pt
\font\eightex=cmex10 at 8pt
\font\eightsy=cmsy9 at 8pt
\font\eightbf=cmbx9 at 8pt
\font\eighttt=cmtt9 at 8pt
\font\eightsl=cmsl9 at 8pt
\font\eightit=cmti9 at 8pt
\font\sixrm=cmr7
\font\sixi=cmmi7
\font\sixex=cmex10 at 7pt
\font\sixsy=cmsy7
\font\sixbf=cmbx7
\font\sixtt=cmtt9 at 7pt
\font\sixsl=cmsl8 at 7pt
\font\sixit=cmti7
\def\rm{\fam0\twelverm}%
  \textfont0=\twelverm \scriptfont0=\eightrm \scriptscriptfont0=\sixrm
  \textfont1=\twelvei \scriptfont1=\eighti \scriptscriptfont1=\sixi
  \textfont2=\twelvesy \scriptfont2=\eightsy \scriptscriptfont2=\sixsy
  \textfont3=\twelveex \scriptfont3=\eightex \scriptscriptfont3=\sixex
  \def\it{\fam\itfam\twelveit}%
  \textfont\itfam=\twelveit
  \def\sl{\fam\slfam\twelvesl}%
  \textfont\slfam=\twelvesl
  \def\bf{\fam\bffam\twelvebf}%
  \textfont\bffam=\twelvebf \scriptfont\bffam=\eightbf
   \scriptscriptfont\bffam=\sixbf
  \def\mit{\fam1}
  \def\oldstyle{\fam1\twelvei}
  \def\tt{\fam\ttfam\twelvett}%
  \textfont\ttfam=\twelvett
  \tt \ttglue=.5em plus.25em minus.15em
\hyphenchar\twelvett=-1 
  \baselineskip=12pt
  \normalbaselineskip=12pt
  \setbox\strutbox=\hbox{\vrule height8.5pt depth3.5pt width\z@}%
  \normalbaselines\rm}

\newdimen\pagewidth \newdimen\pageheight \newdimen\ruleht
 \hsize=38pc  \vsize=51pc  \maxdepth=2.2pt  \parindent=3pc
\newdimen\HSIZE
\HSIZE=\hsize  
\pagewidth=\hsize \pageheight=\vsize \ruleht=.5pt
\abovedisplayskip=6pt plus 3pt minus 1pt
\belowdisplayskip=6pt plus 3pt minus 1pt
\abovedisplayshortskip=0pt plus 3pt
\belowdisplayshortskip=4pt plus 3pt
\overfullrule=0in  
\newlinechar = `\|%
\outer\def\section#1\par{\vskip0pt plus.3\vsize\penalty-23
  \vskip0pt plus-.3\vsize\bigskip\medskip
  \message{|#1|}\leftline{\bf#1}\nobreak\medskip}
\outer\def\subsection#1\par{\vskip0pt plus.3\vsize\penalty-18
  \vskip0pt plus-.3\vsize\bigskip
  \message{|#1|}\noindent{\bf#1}}
\outer\def\proclaim #1. #2\par{\vskip0pt plus.3\vsize\penalty-15
  \vskip0pt plus-.3\vsize\bigskip\vskip\parskip
  \noindent{\sc#1.\enspace}{\sl#2}\par
  \ifdim\lastskip<\medskipamount \removelastskip\penalty13\bigskip\fi}
\outer\def\intermezzo#1{\vskip0pt plus.3\vsize\penalty-10
  \vskip0pt plus-.3\vsize\bigskip\vskip\parskip
  \noindent{\sc#1\enspace}}
\let\im\intermezzo
\outer\def\proof{\vskip0pt plus.3\vsize\penalty-15
  \vskip0pt plus-.3\vsize\bigskip\vskip\parskip
  \noindent{\sc Proof.\enspace}}
\outer\def\remark #1. #2\par{\vskip0pt plus.3\vsize\penalty-10
  \vskip0pt plus-.3\vsize\bigskip\medskip\vskip\parskip
  \noindent{\sc Remark #1.\enspace} #2\par
  \ifdim\lastskip<\medskipamount \removelastskip\penalty13\medskip\fi}
\def\title#1#2#3{\begingroup\null\bigskip\bigskip\tifont\halign{\centerline{##}\cr
  #1\cr #2\cr #3\cr}\bigskip\endgroup}

\outer\def\proclaim #1. #2\par{\vskip0pt plus.1\vsize\penalty-15
  \vskip0pt plus-.1\vsize\bigskip\vskip\parskip
  \noindent{\sc#1.\enspace}{\sl#2}\par
  \ifdim\lastskip<\medskipamount \removelastskip\penalty5\bigskip\fi}
\outer\def\intermezzo#1{\vskip0pt plus.1\vsize\penalty-5
  \vskip0pt plus-.1\vsize\bigskip\vskip\parskip
  \noindent{\sc#1\enspace}}
\let\im\intermezzo
\outer\def\proof{\vskip0pt plus.1\vsize\penalty-10
  \vskip0pt plus-.1\vsize\bigskip\vskip\parskip
  \noindent{\sc Proof.\enspace}}
\outer\def\remark #1. #2\par{\vskip0pt plus.1\vsize\penalty-5
  \vskip0pt plus-.1\vsize\bigskip\medskip\vskip\parskip
  \noindent{\sc Remark #1.\enspace} #2\par
  \ifdim\lastskip<\medskipamount \removelastskip\penalty5\medskip\fi}
\outer\def\abstract{\vskip0pt plus.1\vsize\penalty-25
  \vskip0pt plus-.1\vsize\bigskip\vskip\parskip
  \noindent{\sc Abstract.\enspace}}
\outer\def\continue{\medskip\noindent}
\outer\def\begindisplay{\begingroup\vskip0pt plus.1\vsize\penalty-25
  \vskip0pt plus-.1\vsize\bigskip\vskip\parskip
  \noindent\obeylines}
\outer\def\enddisplay{\ifdim\lastskip<\medskipamount \removelastskip\penalty5\bigskip\fi
  \endgroup}
\let\begintable=\beginalgo
\let\endtable=\endalgo
\def\title#1#2#3{\begingroup\null\bigskip\bigskip\tifont\halign{\centerline{##}\cr
  #1\cr #2\cr #3\cr}\bigskip\endgroup}
\def\devroye{\halign{\hfil##\hfil\cr
				    Luc Devroye\cr
                                     \cr
                                    School of Computer Science\cr
                                    McGill University\cr
                                    3450 University Street\cr
                                    Montreal H3A 2K6\cr
                                    Canada\cr
                                    luc@cs.mcgill.ca\cr
				    \cr}}
\def\janson{\halign{\hfil##\hfil\cr
				    Svante Janson\cr
                                     \cr
				    Matematiska Institutionen \cr
				    Uppsala Universitet \cr
				    Box 480, 751 06 Uppsala \cr
				    Sweden \cr
				    Svante.Janson@math.uu.se \cr
				    \cr }}
\def\author{
	\centerline{\vtop{\devroye}   \phantom{and}   \vtop{\janson}}
	\centerline{}
	\centerline{May 28, 2009}
	\centerline{}
	\centerline{}
}
\outer\def\keywords{\medskip\noindent{\sc Keywords and phrases.}\ }
\outer\def\ams{\medskip\noindent{\sc 2000 Mathematics Subject Classifications:\ }}
\outer\def\crcat{\medskip\noindent{\sc CR Categories:\ }}
\outer\def\endpage{{\nopagenumbers\vfill
  \hrule\smallskip\noindent\ninerm
  \normalbaselineskip=11pt
  \setbox\strutbox=\hbox{\vrule height7.9pt depth3.1pt width\z@}
  \normalbaselines
  The first author's research was sponsored by NSERC Grant A3456.
  The research was mostly done at the Institute Mittag-Leffler during the
  programme ``Discrete Probability'' held in 2009.
  \eject}}
\outer\def\doublespacing{}
\outer\def\singlespacing{
  \normalbaselineskip=13pt
  \setbox\strutbox=\hbox{\vrule height9pt depth4pt width\z@}
  \normalbaselines}

\def\beginalgo{\begingroup\vskip0pt plus .1\vsize\penalty-18
  \vskip0pt plus-.1\vsize\medskip\vskip\parskip\tt}
\def\endalgo{\par\bigskip
  \ifdim\lastskip<\medskipamount \removelastskip\penalty5\medskip\fi
  \endgroup}

\twelvepoint  
\parskip=14pt
\baselineskip=14pt
\normalbaselineskip=14pt
\setbox\strutbox=\hbox{\vrule height10pt depth4pt width\z@}%
\normalbaselines
\raggedbottom

\title{Long and short paths in uniform random recursive dags}{}{}

$$\hbox to \hsize{\hfill\vbox{
\author}\hfill
}$$
  \bigskip

\vfill

\abstract
In a uniform random recursive $k$-dag, there is a
root, $0$, and each node in turn,
from $1$ to $n$, chooses $k$ uniform random parents from
among the nodes of smaller index. 
If $S_n$  is the shortest path  distance from node $n$
to the root, then we determine the constant $\sigma$ such
that $S_n / \log n \to \sigma$ in probability as $n \to \infty$.
We also show that $\max_{1 \le i \le n} S_i / \log n \to \sigma$
in probability.

\vfill

\keywords
Uniform random recursive dag.
Randomly generated circuit.
Random web model.
Longest paths.
Probabilistic analysis of algorithms.
Branching process.
\vfill
\crcat
F.2.0, G.3.0, G.2.2
\vfill
\ams
60C05, 60F99, 68R01.

\vfill
\endpage
\smallskip

\vfill\eject

\section
1. Introduction.

A uniform random $k$-dag is an infinite directed graph
defined as follows.
For each  of the integers $1, 2, \ldots$, 
we pick a random set of $k$ parents with replacement
uniformly from among the smaller non-negative integers.
This defines an infinite directed acyclic graph (or, dag)
with one root (0),
and can be viewed as a (too) simplistic model of the web,
a random recursive circuit (Diaz, Sperna, Spirakis, Toran and Tsukiji, 1994,
and Tsukiji and Xhafa, 1996),
and a generalization of the {\sc urrt} (uniform random recursive tree),
which is obtained for $k=1$.
All the asymptotic results in the paper remain valid
when parents are selected without replacement.

The uniform random $k$-dag restricted to vertices
$0, 1, \ldots, n$, is denoted by $U_{k,n}$ or simply
$U_n$. Indeed, we will take $k=2$ in the main part of the paper,
and point out the obvious modifications needed when $k>2$ as we proceed.
The infinite dag is denoted by $U_\infty$.

From a given node $n$, let ${\cal P}_n$
be the collection of paths from node $n$ to the origin.
The length of path $p \in {\cal P}_n$ is $L(p)$.
One can consider various path lengths:
$$
S_n = \min_{p \in {\cal P}_n} L(p) ~,~
R^-_n = L(P^-_n) ~,~
R_n = L(P_n) ~,~
R^+_n = L(P^+_n) ~,~
L_n = \max_{p \in {\cal P}_n} L(p),
$$
where $S$, $R$ and $L$ are mnemonics for shortest, random, and longest,
and $P^-_n, P_n$ and $P^+_n$ are the paths in ${\cal P}_n$, where we follow
the parent with the smallest index, the first parent and the parent with the largest index,
respectively.
We can regard $R^-_n$ and $R^+_n$ as greedy
approximations of $S_n$ and $L_n$ respectively. Note that,
at least in a stochastic sense,
$$
S_n \le R^-_n \le R_n \le R^+_n \le L_n.
$$

The length of the longest path is relevant for the time
to compute the value of node $n$ in a random recursive circuit,
when nodes know their value only when all parents know their value.
However, there are situations in which node values are determined
as soon as one parent or a subset of parents know their value---they
are called self-time circuits by Codenotti, Gemmell and Simon (1995).
For the one-parent case, this leads naturally to the study of $S_n$.
In networks, in general, shortest paths have been of interest
almost since they were conceived (Prim, 1957; Dijkstra, 1959).

It is of interest to study the extreme behavior, as measured by
$$
\max_{1 \le \ell \le n} S_\ell ~,~ 
\max_{1 \le \ell \le n} R^-_\ell ~,~ 
\max_{1 \le \ell \le n} R_\ell ~,~ 
\max_{1 \le \ell \le n} R^+_\ell ~,~ 
\max_{1 \le \ell \le n} L_\ell.
$$
If we replace max by min in these definitions, we obtain the constant $1$,
and it is therefore more meaningful to ask for the exteme minimal
behavior as defined by
$$
\min_{n/2 \le \ell \le n} S_\ell ~,~ 
\min_{n/2 \le \ell \le n} R^-_\ell ~,~ 
\min_{n/2 \le \ell \le n} R_\ell ~,~ 
\min_{n/2 \le \ell \le n} R^+_\ell ~,~ 
\min_{n/2 \le \ell \le n} L_\ell.
$$
So, in all, there are fifteen parameters that could be studied.

We take this opportunity to introduce the label process,
which
will be referred to throughout the paper.
The label of each parent of $n$ is distributed
as $\lfloor nU \rfloor$, with $U$ uniform $[0,1]$.
An $\ell$-th generation ancestor has a label distributed like
$$
\lfloor \cdots \lfloor \lfloor n U_1 \rfloor U_2 \rfloor \cdots U_\ell \rfloor
\in
\left[
n U_1 U_2 \cdots U_\ell - \ell ,
n U_1 U_2 \cdots U_\ell
\right] ,
$$
where the $U_i$'s are i.i.d.\ uniform $[0,1]$ random variables.

\noindent
{\sc The parameter $R_n$.}
It is clear that $R_n$ is just the distance from node $n$
in a {\sc urrt} to its root.
In particular, $R_n$ and its minimal and maximal versions do
not depend upon $k$.
We dispense immediately with $R_n$ and its extensions because
of well-known results on the {\sc urrt} obtained via the
study of branching random walks by Devroye (1987)
and the equivalence between $R_n$ and the number of records
in an i.i.d.\ sequence of continuous random variables (see, e.g.,
R\'enyi (1962), Pyke (1965), Glick (1975) or Devroye (1988)).
Only the minimal parameter for $R_n$ requires a gentle
intervention.
We know that
$$
{R_n \over \log n} \to 1 ~\hbox{\rm in probability},
$$
for example. Furthermore, 
$$
{R_n - \log n \over \sqrt{\log n}} \tendsinlaw {\cal N} 
$$
where $\cal N$ is a standard normal random variable, 
and $\tendsinlaw$ denotes convergence in distribution.
Furthermore, an explicit tail bound on $R_n$ will be needed further
on in the paper.
The maximal value of $R_\ell, 1 \le \ell \le n$,
follows immediately from either Devroye (1987) or Pittel (1994).
We summarize:

\proclaim Theorem 1.
We have
$$
{R_n \over \log n} \to 1 ~\hbox{\rm in probability},
$$
$$
{\max_{1 \le \ell \le n} R_\ell \over \log n} \to e ~\hbox{\rm in probability},
$$
and
$$
\lim_{n \to \infty} \PROB \left\{ \min_{n/2 \le \ell \le n} R_\ell \le 2 \right\}
= 1.
$$
Finally, for $t \ge \log n $ integer,
$$
\PROB \{ R_n > t \}
\le \exp \left( t - \log n -t \log (t/\log n) \right).
$$

\proof
An outline of proof is needed for the third part
and the explicit bound in part four.
Let us count the number of nodes of index in $[1,n/2]$ that
connect directly to the root. This number is
$$
Z = \sum_{\ell = 1 }^{n/2} \xi_{1/\ell},
$$
where $\xi_p$ is Bernoulli $(p)$. 
Let $A$ be the event that no node of index in $(n/2, n]$
connects to a node counted in $Z$. This probability
is smaller than
$$
\eqalignno{
\EXP \left\{ \left( 1 - {Z \over n} \right)^{n/2} \right\}
&\le \EXP \left\{ e^{-Z/2} \right\} \cr
&\le \prod_{\ell=1}^{n/2} \left( 1 - 1/\ell + 1/(\sqrt{e}\ell) \right)  \cr
&\le \exp \left( - \sum_{\ell=1}^{n/2} {1-1/\sqrt{e} \over \ell} \right)  \cr
&\le \left( \lfloor n/2 \rfloor \right)^{1-1/\sqrt{e}}. \cr}
$$
If the complement of $A$ holds, then clearly, $\min_{n/2 \le \ell \le n} R_\ell \le 2$,
and thus, we have shown the third part of Theorem 1.
Turning to part four,
note that $R_n \le \min \{ t : n U_1 \cdots U_t < 1 \}$,
and thus that
$$
\PROB \{ R_n > t \}
\le \inf_{\lambda > 0} \EXP \left\{ (n U_1 \cdots U_t )^\lambda \right\}
= \inf_{\lambda > 0} n^\lambda (\lambda+1)^{-t}
= \exp \left( t - \log n -t \log (t/\log n) \right). ~\square
$$
\medskip

\proclaim Conjecture.
For all fifteen parameters, generically denoted by
$X_n$, there exist finite constants $x = x(k) \ge 0$ such that
$$
{X_n \over \log n} \to x
~ \hbox{\rm in probability}.
$$

\noindent {\sc Remark.}
The limits in the conjecture are denoted by $\sigma$, $\rho^-$, $\rho$, $\rho^+$
and $\lambda$ for $S_n$, $R^-_n$, $R_n$, $R^+_n$ and $L_n$,
respectively. For the minimal and maximal versions of these
parameters, we will use the subscripts $\min$ and $\max$,
respectively, as in $\rho^+_{\min{}}$ and $\sigma_{\max{}}$, for example.

%
\bigskip

Let us briefly survey what is known and provide conjectures
in the other cases.

\noindent
{\sc The parameter $L_n$.}
Tsukiji and Xhafa (1996) showed that $\lambda_{\max{}} = ke$.
The Chernoff large deviation bound shows that
$\lambda$ is at most the largest solution $x$ of
$$
\left( {ke \over x } \right)^x e^{-1} = 1,
\hfil
\eqno{(1)}
$$
and thus $\lambda < \lambda_{\max{}}$.
We believe that $\lambda$ is indeed given by (1) based
on arguments not unlike the proof of Theorem 2 below.
We have no guess at this point about the value of $\lambda_{\min{}}$.
\bigskip

\noindent
{\sc The parameter $R^+_n$.}
In the label process,
the parent's index is approximately 
distributed as $n \max (U_1, \ldots , U_k)$, where the
$U_i$'s are i.i.d.\ uniform $[0,1]$ random variables.
If $U$, as elsewhere in this paper, is uniform $[0,1]$,
then the parent's index is thus roughly like $n U^{1/k}$.
By renewal theory, this implies that 
$$
{R^+_n \over \log n} \to k \isdef \rho^+ ~\hbox{\rm in probability}.
$$
Chernoff's large deviation bound show that $\rho^+_{\max{}}$
is at most the unique solution $x$ of (2) that is above $k$:
$$
\left( {ke \over x } \right)^x e^{1-k} = 1.
\hfil
\eqno{(2)}
$$
We believe that the solution of (2) yields $\rho^+_{\max{}}$.
Applying Chernoff to the other tail shows that
$\rho^+_{\min{}}$ is at least the other solution of (2),
as (2) has two solutions, one below $k$ and one above $k$.
Furthermore, we believe that the two solutions  of (2)
yield $\rho^+_{\min{}}$ and $\rho^+_{\max{}}$.

For $k=2$, the parameter $R^+_n$  is intimately linked to
the random binary search tree, which can be grown
incrementally by a well-known process described as follows: given
an $n$-node random binary search tree, sample one of its $n+1$ 
external nodes uniformly at random, replace it by node $n+1$, and
continue. The parent of that node is either its neighbor (in the total
ordering) to the left or its neighbor to the right, and in fact, it is
the neighbor added last to the tree. But
the labels (times of insertion) of the neighbors are uniformly drawn
without replacement from $\{1,\ldots,n\}$, and are thus
roughly distributed as $n U$, so that the parent of $n+1$
is roughly distributed as $n \sqrt{U}$, because the maximum of
two i.i.d.\ uniform $[0,1]$ random variables is distributed as $\sqrt{U}$.
With this in mind, $\max_{1 \le \ell \le n} R^+_\ell$ is the height
of the random binary search tree, $R^+_n$ is the depth (distance to the root)
of the node of label $n$ 
(the $n$-th node inserted), and $\min_{n/2 \le \ell \le n} R^+_\ell$
is very roughly the shortest distance from leaf to root, or fill-up level.
These quantities behave in probability as described above,
as shown by Devroye (1986, 1987), and this
explains the values $\rho^+_{\max{}} = 4.31107\ldots$,
$\rho^+ = 2$ and $\rho^+_{\min{}} = 0.3733\ldots$.
\bigskip

\noindent
{\sc The parameter $R^-_n$.}
Arguing as above, the parent's index is approximately
distributed as $n \min (U_1,\ldots,U_k)$.
By a property of the uniform (or exponential) distribution,
using a sequence of i.i.d.\ exponential random variables $E_1,E_2, \ldots$,
we have this distributional identity:
$$
n \min (U_1,\ldots,U_k)
\inlaw n U_1 U_2^{1/2} \cdots U_k^{1/k} 
\inlaw \exp \left( \log n - \sum_{j=1}^k {E_j \over j} \right).
$$
Renewal theory easily gives the law of large numbers and
central limit theorem for $R^-_n$. For example,
$$
{R^-_n \over \log n} \to {1 \over H_k} \isdef \rho^- ~\hbox{\rm in probability},
$$
where $H_k = \sum_{j=1}^k (1/j)$ is the $k$-th harmonic number.
Using large deviation bounds similar to the ones
used below in showing part of Theorem 2, one gets that
$$
\lim_{n \to \infty} \PROB \left\{
\max_{1 \le \ell \le n} R^-_\ell \ge (x + \epsilon) \log n \right\}
= 0
$$
for all $\epsilon > 0$,
where $x$ is the solution greater than $1/H_k$ of
$$
1 + f(x) = x \sum_{j=1}^k \log \left( 1 + f(x)/j \right),
$$
and $f(x) > 0$ is implicitly defined by
$$
\sum_{j=1}^k { 1 \over j+f(x) } = {1 \over x }, \qquad x > 1/H_k.
$$
These equations follow from the obvious Chernoff bound.
We conjecture that
$\rho^-_{\max{}}$ equals this upper bound, but a rigorous proof that $\rho^-_{\max{}}$
is indeed as described above is not given in this paper.
\bigskip

\noindent
{\sc The parameter $S_n$.}
The most important parameter for computer scientists and combinatorialists
is the one in which graph distances are defined by shortest paths, and this
leads to the study of $S_n$. That was the original motivation of the paper,
and we will settle first order asymptotics in this paper.
Theorem 1 implies, for example, that with probability tending to one,
$$
\min_{n/2 \le \ell \le n} S_\ell \le 2.
$$
So we turn to $\sigma$ and $\sigma_{\max{}}$:

\proclaim Theorem 2.
Assume $k \ge 2$.
Then $\sigma = \sigma_{\max{}}$, where $\sigma$ is given by the solution $x \in (0,1)$ of
$$
\varphi(x) \isdef \left( { ke \over x} \right)^x  e^{-1} = 1.
\hfil
\eqno{(3)}
$$
[Note that $\varphi$ is indeed an increasing function on $(0,1)$.]

Observe that Theorem 2 does not extend to $k=1$, because in that case,
$S_n \equiv R_n \equiv L_n$, and similarly for the maximal versions of these
parameters, in view of the equivalence with the {\sc urrt}.
Thus, $S_n / \log n \to 1$ and $\max_{1 \le \ell \le n} S_\ell / \log n \to e$
in probability.
\bigskip

The following is a table of constants in the Conjecture for
$k=2$. The constants involving $\sigma$ (top row) are obtained
in this paper, while those involving $\rho$ (third row) are covered by
Theorem 1. The constants $\rho^-$ and $\rho^+$ follow from
ordinary renewal theory. The zeroes in the table follow
from Theorem 1. Finally, $\lambda_{\max{}}$ is due
to  Tsukiji and Xhafa (1996). There are thus four conjectured constants,
which happen to be one-sided bounds ($\rho^-_{\max{}}$, $\rho^+_{\min{}}$,
$\rho^+_{\max{}}$, $\lambda$),  and one unknown constant, $\lambda_{\min{}}$.

\bigskip
$$\vbox{
\halign{\quad#\quad&\quad#\hfil\quad&\quad#\hfil\quad&\quad#\hfil\quad\cr
\noalign{\smallskip}
\noalign{\hrule}
\noalign{\smallskip}
$\sigma_{\min{}} \hfil \sigma \hfil \sigma_{\max{}}$&
0 & 0.3733$\ldots$  & 0.3733$\ldots$ \cr
$\rho^-_{\min{}} \hfil \rho^- \hfil \rho^-_{\max{}}$&
0 & 0.6666$\ldots$ ($= 2/3$) &  1.6737$\ldots$ \cr
$\rho_{\min{}} \hfil \rho \hfil \rho_{\max{}}$&
0 & $1$ &  2.7182$\ldots$ ($=e$)\cr
$\rho^+_{\min{}} \hfil \rho^+ \hfil \rho^+_{\max{}}$&
0.3733$\ldots$ & $2$ &  4.3110$\ldots$ \cr
$\lambda_{\min{}} \hfil \lambda \hfil \lambda_{\max{}}$&
? & 4.3110$\ldots$ &  5.4365$\ldots$ ($=2e$)\cr
\noalign{\smallskip}
\noalign{\hrule}
\noalign{\smallskip}
}}$$
\bigskip

\medskip

The following is a table of $\sigma$, $\rho^-$ and $\rho^-_{\max{}}$
for different numbers of parents $k$.

\bigskip
$$\vbox{
\halign{\hfil#\quad&\quad#\hfil\quad&\quad#\hfil\quad&\quad#\hfil\quad\cr
$k$&\hfil $\sigma$&\hfil $\rho^-$&\hfil $\rho^-_{\max{}}$\cr
\noalign{\smallskip}
\noalign{\hrule}
\noalign{\smallskip}
2&
0.3733$\ldots$&
0.6666$\ldots$&
1.6737$\ldots$
\cr
3&
0.3040$\ldots$&
0.5454$\ldots$&
1.3025$\ldots$
\cr
4&
0.2708$\ldots$&
0.48&
1.1060$\ldots$
\cr
5&
0.2503$\ldots$&
0.4379$\ldots$&
0.9818$\ldots$
\cr
6&
0.2361$\ldots$&
0.4081$\ldots$&
0.8951$\ldots$
\cr
7&
0.2254$\ldots$&
0.3856$\ldots$&
0.8305$\ldots$
\cr
8&
0.2170$\ldots$&
0.3679$\ldots$&
0.7800$\ldots$
\cr
9&
0.2102$\ldots$&
0.3534$\ldots$&
0.7393$\ldots$
\cr
10&
0.2045$\ldots$&
0.3414$\ldots$&
0.7057$\ldots$
\cr
11&
0.1996$\ldots$&
0.3311$\ldots$&
0.6773$\ldots$
\cr
12&
0.1954$\ldots$&
0.3222$\ldots$&
0.6531$\ldots$
\cr
13&
0.1916$\ldots$&
0.3144$\ldots$&
0.6318$\ldots$
\cr
14&
0.1883$\ldots$&
0.3075$\ldots$&
0.6132$\ldots$
\cr
15&
0.1854$\ldots$&
0.3013$\ldots$&
0.5966$\ldots$
\cr
16&
0.1827$\ldots$&
0.2957$\ldots$&
0.5816$\ldots$
\cr
17&
0.1802$\ldots$&
0.2907$\ldots$&
0.5683$\ldots$
\cr
18&
0.1780$\ldots$&
0.2861$\ldots$&
0.5560$\ldots$
\cr
19&
0.1760$\ldots$&
0.2818$\ldots$&
0.5448$\ldots$
\cr
20&
0.1740$\ldots$&
0.2779$\ldots$&
0.5346$\ldots$
\cr
21&
0.1723$\ldots$&
0.2743$\ldots$&
0.5251$\ldots$
\cr
22&
0.1706$\ldots$&
0.2709$\ldots$&
0.5164$\ldots$
\cr
23&
0.1691$\ldots$&
0.2677$\ldots$&
0.5083$\ldots$
\cr
24&
0.1676$\ldots$&
0.2648$\ldots$&
0.5007$\ldots$
\cr
25&
0.1663$\ldots$&
0.2620$\ldots$&
0.4936$\ldots$
\cr
26&
0.1650$\ldots$&
0.2594$\ldots$&
0.4868$\ldots$
\cr
27&
0.1638$\ldots$&
0.2569$\ldots$&
0.4805$\ldots$
\cr
28&
0.1626$\ldots$&
0.2546$\ldots$&
0.4747$\ldots$
\cr
29&
0.1615$\ldots$&
0.2524$\ldots$&
0.4690$\ldots$
\cr
30&
0.1604$\ldots$&
0.2503$\ldots$&
0.4638$\ldots$
\cr
35&
0.1559$\ldots$&
0.2411$\ldots$&
0.4409$\ldots$
\cr
40&
0.1521$\ldots$&
0.2337$\ldots$&
0.4225$\ldots$
\cr
45&
0.1490$\ldots$&
0.2275$\ldots$&
0.4074$\ldots$
\cr
50&
0.1463$\ldots$&
0.2222$\ldots$&
0.3946$\ldots$
\cr
\noalign{\smallskip}
\noalign{\hrule}
\noalign{\smallskip}
}}$$
\bigskip

\medskip

\section
2. The shortest path length $S_n$.

We will establish Theorem 2
in two parts. First we show that for all $\epsilon > 0$,
$$
\lim_{n \to\infty} \PROB \{ S_n \le (1-\epsilon) \sigma \log n \} = 0,
\hfil\eqno{(4)}
$$
and then that
$$
\lim_{n \to\infty} \PROB \left\{ \max_{1 \le \ell \le n} S_\ell \ge (1+\epsilon) \sigma \log n \right\} = 0.
\hfil\eqno{(5)}
$$
We only consider the case $k=2$ since the case $k > 2$ follows quite
easily.

\proclaim Lemma 1.
Let $G_a$ be gamma$(a)$, with $a \ge 1$. Then
$$
{ \PROB \{ G_a \ge x \} \over { x^{a-1} e^{-x} \over \Gamma(a) }} \le {1 \over 1 - {a-1 \over x } }  , x > a-1,
$$
and
$$
{ \PROB \{ G_a \le x \} \over { x^{a-1} e^{-x} \over \Gamma(a) }} \le  {1 \over {a-1 \over x } - 1 } , x < a-1.
$$

\proof
The gamma density is $f(y) = y^{a-1} e^{-y}/\Gamma (a)$. It is log-concave for $a\ge 1$, and
thus, a first-term Taylor series bound yields the inequality
$$
f(y) \le f(x) e^{(y-x) (\log f)' (x)} = f(x) e^{(y-x) ( (a-1)/x \,- 1 ) }.
$$
Integrating the upper bound out over $[x,\infty)$ or $(-\infty,x]$ then immediately
yields the results. $\square$
\medskip


From node $n$, we can consider the index of the first of
the $2^\ell$  $\ell$-th level ancestors, which is distributed as
$$
\lfloor \cdots \lfloor \lfloor n U_1 \rfloor U_2 \rfloor \cdots U_\ell \rfloor
\ge
n U_1 U_2 \cdots U_\ell - \ell
\inlaw
n \exp (-G_\ell) - \ell,
$$
where $\inlaw$ denotes equality in distribution, and $G_\ell$ is gamma$(\ell)$.
If these indices are $I_1,\ldots, I_{2^\ell}$, then we have
$$
\eqalignno{
\PROB \{ S_n \le \ell \} 
&= \PROB \left\{ \min_{1 \le i \le 2^\ell} I_i = 0 \right\} \cr
&\le 2^\ell \PROB \{ I_1 = 0 \} \cr
&\le 2^\ell \PROB \{ n \exp (-G_\ell) - \ell \le 0 \} \cr
&= 2^\ell \PROB \{ G_\ell \ge \log (n/\ell) \} \cr
&\le { 2^\ell (\log (n/\ell))^{\ell-1} e^{-\log (n/\ell)} \over \Gamma(\ell) \left( 1 - {\ell-1 \over \log (n/\ell) } \right) } \qquad \hbox{\rm (if $\log (n/\ell) \ge \ell -1$)} \cr
&\le { \ell^{3/2} (2 \log (n))^{\ell} e^{-\log (n)} \over (\ell/e)^\ell \left( 1 - {\ell-1 \over \log (n/\ell) } \right)  } . \cr}
$$
Set $\ell = \lfloor t \log n \rfloor$ for $t \in (0,1)$, and note that the upper bound is
$$
\Theta \left( \log^{3/2} (n) \right)
\times (\varphi (t))^{\log n},
$$
where $\varphi(t) = (2e/t)^t / e$ is as in (3).
We have $\varphi (\sigma) = 1$ for $\sigma = 0.3733\ldots$.
Thus, we have shown (4): for all $\epsilon > 0$,
$$
\PROB \{ S_n \le (\sigma - \epsilon) \log n  \} = o(1).
$$

Although we will not need it directly, we will also deal with
the upper bound on $S_n$.
This can be done in a number of ways, but the shortest route is perhaps via
the great-grandparent strategy that jumps $\ell$ generations at a time,
where $\ell$ now is a large but fixed integer.
We denote this by $\ell$-{\sc ggp}.
We associate with each node $n$ two independent uniform $[0,1]$
integers $U$ and $V$ and let the parent labels be $\lfloor nU \rfloor$
and $\lfloor nV \rfloor$.
Let $A_n$ be the event that any of the $2^\ell$ ancestors of node $n$
conicide. It is clear that $\PROB \{ A_n \} \to 0$ as $n \to \infty$.
As an ancestor label is described by
$$
\lfloor \cdots \lfloor \lfloor n U_1 \rfloor U_2 \rfloor \cdots U_\ell \rfloor
\le
n U_1 U_2 \cdots U_\ell  \inlaw n \exp (-G_\ell),
$$
we define
$$
Z_\ell = \min_{p \in \cal P} \prod_{e \in p} U_e
$$
where ${\cal P}$ is the collection of all paths of length $\ell$
above node $n$, and each $p \in {\cal P}$ consists of edges $e$ that
each have an independent uniform random variable associated with it.
If $\epsilon > 0$ and $n$ is greater than some $n_\epsilon$,
then the $\ell$-{\sc ggp} gives with probability greater than $1-\epsilon$
a node with label less than $Z_\ell n$.
Define
$$
Z_\ell^{(\epsilon)} = \min (Z_\ell, b)
$$
where $b$ is chosen such that $\PROB \{ Z_\ell > b \} = \epsilon$.
As long as the label stays above $n_\epsilon$, one can dominate the
labels in the $\ell$-{\sc ggp} by multiplying $n$ with successive
independent copies of $Z_\ell^{(\epsilon)}$.
Let $T_n$ be the number of steps until the label in $\ell$-{\sc ggp}
reaches $n_\epsilon$ or less.
Renewal theory shows that with probability tending to one,
$$
T_n \le {(1 + \epsilon) \log n \over \EXP \left\{ - \log \left( Z_\ell^{(\epsilon)} \right) \right\} }.
$$
because the $\ell$-{\sc ggp} takes $\ell$ steps at a time,
and because a node with label $n_\epsilon$ is not further than $n_\epsilon$
away from the origin, we see that with probability tending to one,
$$
\eqalignno{
S_n 
&\le n_\epsilon +  {\ell (1 + \epsilon) \log n \over \EXP \left\{ - \log \left( Z_\ell^{(\epsilon)} \right) \right\}} \cr
&\le {\ell (1 + 2\epsilon) \log n \over \EXP \left\{ - \log \left( Z_\ell^{(\epsilon)} \right) \right\} }.
\cr}
$$

Uniform integrability implies that
$$
\lim_{\epsilon \downarrow 0} \EXP \left\{ - \log \left( Z_\ell^{(\epsilon)} \right) \right\}
= \EXP \left\{ - \log \left( Z_\ell \right) \right\} .
$$
Therefore, for any (new, fresh) $\epsilon > 0$ and $\ell \ge 1$,
with probability going to one,
$$
S_n  \le {\ell (1 + \epsilon) \log n \over \EXP \left\{ - \log \left( Z_\ell \right)  \right\} }.
$$
Observe that
$$
{ - \log \left( Z_\ell \right) \over \ell }
\inlaw
{1 \over \ell}  \max_{p \in {\cal P}} \sum_{e \in p} E_u,
$$
where the $E_u$ are i.i.d.\  exponential random variables.
From the theory of branching random walks, it is easy to
verify (see, e.g., Biggins (1977), or Devroye (1986, 1987))
that, as $\ell \to \infty$,
$$
{1 \over \ell}  \max_{p \in {\cal P}} \sum_{e \in p} E_u 
\to {1 \over \sigma}
$$
in probability. Thus,
$$
\liminf_{\ell \to \infty} { - \EXP \left\{ \log \left( Z_\ell\right) \right\} \over \ell } \ge { 1 \over \sigma },
$$
and thus, by choosing $\ell$ large enough, we see that with probability tending to one,
$$
S_n  \le  (1 + 2 \epsilon) \sigma \log n .
$$
This concludes the proof of the first part of Theorem 2.

The next section requires an explicit rate of convergence.
To this end, still restricting ourselves to $k=2$ only,
let $Z_{\ell,1}^{(\epsilon)}, Z_{\ell,2}^{(\epsilon)}, \ldots$ be i.i.d.\ copies of 
$Z_\ell^{(\epsilon)}$, and
note that, 
$$
\eqalignno{
T_n 
&\le \min \left\{ t :  n Z_{\ell,1}^{(\epsilon)} \cdots Z_{\ell,t}^{(\epsilon)} < 1 \right\} \cr
&= \min \left\{ t :   \log \left( 1/Z_{\ell, 1}^{(\epsilon)} \right) + \cdots + \log \left( 1/Z_{\ell, t}^{(\epsilon)} \right) > \log n \right\}.
\cr}
$$
Set $\mu =  \EXP \left\{ \log \left( 1/Z_\ell^{(\epsilon)} \right) \right\}$.
Then, assuming $\delta^* \in (0, 1/2)$ and $\delta \in (\delta^* , 2 \delta^*)$
such that $m = (1/\mu + \delta ) \log n$ is integer-valued,
$$
\eqalignno{
\PROB \{ T_n >  m \}
&\le \PROB \left\{ \log \left( 1/Z_{\ell,1}^{(\epsilon)} \right) + \cdots + \log \left( 1/Z_{\ell, m}^{(\epsilon)} \right)  < \log n \right\} \cr
&= \PROB \left\{ \log \left( 1/Z_{\ell,1}^{(\epsilon)} \right) + \cdots + \log \left( 1/Z_{\ell, m}^{(\epsilon)} \right)  -m \mu < - \delta \mu \log n \right\}.
\cr}
$$

Let $p > 2$ be a fixed number.
Rosenthal's inequality (Rosenthal, 1970, Fuk and Nagaev, 1971, see also Petrov, 1975)
states that there is a constant $C_p$ with the following property.
If $ \{X_n , n \ge 1 \}$ is  a sequence of centered and independent random 
variables, and if $Y_n = X_1 + \cdots + X_n$, and if $\EXP \{ |X_n |^p \} < \infty$
for all $n$, then
$$
\EXP \{ |Y_n |^p \}
\le
C_p \left(  \sum_{j=1}^n \EXP \{ |X_j |^p \}  + \left( \Var \{ Y_n \} \right)^{p/2} \right).
$$
For i.i.d.\ random variables with $X_1 = X$, we have
$$
\EXP \{ |Y_n |^p \}
\le C_p \left(  n \EXP \{ |X|^p \}  + n^{p/2} \left( \EXP \{ X^2 \}  \right)^{p/2} \right)
\le 2 C_p \max (n, n^{p/2} ) \EXP \{ |X|^p \} .
$$
Applied to our situation with $p=4$, using Markov's inequality, we have
$$
\eqalignno{
\PROB \{ T_n >  m \}
&\le (\delta \mu \log n)^{-4} \EXP  \left\{ \left( \log \left( 1/Z_{\ell,1}^{(\epsilon)} \right) + \cdots + \log \left( 1/Z_{\ell,m}^{(\epsilon)} \right)  -m \right)^4 \right\} \cr
&\le 2 C_4 (\delta \mu \log n)^{-4} m^2 \EXP \left\{ \left|\log \left( 1/Z_\ell^{(\epsilon)} \right) - \mu \right|^4 \right\} \cr
&\le C (\log n)^{-2} {\delta^*}^{-4} ,
\cr}
$$
where $C$ depends upon $\epsilon$ and $\ell$ only.
The remainder of the argument involving an appropriate choice of  $\ell$
remains valid, and we can conclude that for any $\epsilon > 0$,
$$
\PROB \{ S_n > (\sigma + \epsilon) \log n \} = O \left( 1/ \log^2 n \right),
\hfil \eqno{(6)}
$$
with room to spare.

\medskip

\section
3. The maximal shortest path length

The purpose of this section is to show (5).
We let $\sigma$ be as in the first part of the proof, and let $\epsilon > 0$
be arbitrary. 
Fix $n$ large enough.
From (6),
$$
\EXP \left\{ \left| \{ j: n/2 \le j \le n , S_j > (\sigma + \epsilon) \log n \} \right| \right\}
= O \left( { n \over \log^2 n } \right),
$$
and thus $\PROB \{ A(n) \} = O \left( { 1 \over \log^2 n } \right)$, where
$$
A(n) \isdef
\left[ \left| \left\{ j: n/2 \le j \le n , S_j > (\sigma + \epsilon) \log n \} \right|  > {n \over 4} \right\} \right].
$$
If we take an incremental view of the process of adding edges,
then a node with index 
in $[n, 2n]$ selects a parent of depth $\le (\sigma + \epsilon) \log n$
and index $\ge n/2$ with probability $\ge 1/8$
if $A(n)$ fails to hold.
It is this observation that will allow us to uniformly bound all depths
by something close to $(\sigma + \epsilon) \log n$.

Consider the indices in dyadic groups, $\{ 2^{r-1} + 1 , \ldots, 2^r \}$, $r \ge 1$.
We recall from a comparison with the {\sc urrt}, that $S_n \le R_n$
and thus that $\max_{1 \le j \le n} S_j \le \max_{1 \le j \le n} R_j$,
and that  (see Theorem 1)
$$
\PROB \left\{ \max_{1 \le j \le n} R_j > 2e \log n \right\} 
\le n^{-2e \log (2)} < n^{-3}.
$$
Thus, for $\gamma > 0$ small enough,
$$
\PROB \left\{ \max_{1 \le j \le \lfloor n^{\gamma} \rfloor } S_j > (\sigma + \epsilon) \log n \right\} = O(n^{-3\gamma}) = o(1).
$$
It remains to show that
$$
\PROB \left\{ \max_{n^\gamma \le j \le n } S_j > (\sigma + \epsilon) \log n \right\} = o(1).
$$
Consider the event 
$$
B(r) = \bigcup_{r' \le s \le r} A(2^s),
$$
where $r'$ is the largest integer such that $2^{r'} < n^\gamma$.
Clearly, $\PROB \{ B(r) \} = O(1/r') = O(1/ \log n)$.
On the complement, $(B(r))^c$, 
intersected with 
$\left[ \max_{1 \le j \le \lfloor n^{\gamma} \rfloor } S_j \le (\sigma + \epsilon) \log (n) \right]$,
 we look at the process started at a node
$m \le n$ and assume that its index 
$m$ is in $ \left\{ 2^{r} + 1 , \ldots, 2^{r+1} \right\}$.
That process is looked at as a binary tree of consecutive parents,
and will be cut off at height $h = \lfloor 10 \log \log n \rfloor$.
There may be duplicate parents (in which case the tree degenerates
to a dag), so we need to be a bit careful.
If any parent in the tree
is selected with index $\le 2^{r'} < n^\gamma$,
then $S_m \le (\sigma + \epsilon) \log n + h$, 
and thus, we can assume that in this ``tree'' any node $j$ selects its parent uniformly
in the range $(2^{r'}, j)$.
At any stage, by our assumption, the probability of picking a parent $i$
having $S_i \le (\sigma + \epsilon) \log n$ is at least $1/8$  (and this
is why we needed the dyadic trick, so that we can make this statement
regardless of the choice of $i$ within the range $(2^{r'}, n]$).
We claim that this ``tree'' has at least $2^{h-1}$ leaves or reaches
$[1,2^{r'}]$ with overwhelming probability.
To see this, note that a node $j$ in it picks a node already selected
with probability not exceeding $2^h /j$.
But the index $j$ is 
stochastically larger than 
$$
X_h \isdef \lfloor \cdots \lfloor \lfloor m U_1 \rfloor U_2 \rfloor \cdots U_h \rfloor
$$
by our remarks about the labeling process.
The probability that there are in fact at least two such unwanted
parent selections (but none of them less than $n^\gamma$) in that ``tree'' is not more than
$$
2^{2h+2} \times \EXP^2 \left\{ { 2^h  \over  X_h } \IND{X_h \ge n^\gamma} \right\}
\le 2^{4h+2} \times \EXP^2 \left\{ { 1  \over  X_h } \IND{X_h \ge n^\gamma}\right\}
\hfil \eqno{(7)}
$$
We have
$$
\eqalignno{
\EXP \left\{ X_h^{-1}  \IND{X_h \ge n^\gamma} \right\}
&= \int_0^\infty \PROB \{ X_h^{-1}  \IND{X_h \ge n^\gamma} > t \} \, dt \cr
&= \int_0^{1/n^\gamma} \PROB \{ X_h < 1/t  \} \, dt \cr
&\le \int_0^{1/n^\gamma} \PROB \{ m U_1 \cdots U_h  < h + 1/t  \} \, dt \cr
&= \int_0^{1/n^\gamma} \PROB \{ \log {m \over h+1/t}  < G_h \} \, dt  \cr
&= \int_0^{1/n^\gamma} \int_{\log_+ { m \over h+1/t}}^\infty {y^{h-1} e^{-y} \over \Gamma(h) } \, dy \, dt \cr
&= \int_0^\infty {y^{h-1} e^{-y} \over \Gamma(h)}  \min \left( n^{-\gamma} , {1 \over \left( me^{-y} - h \right)_+ } \right)   \, dy \cr
&\le \int_0^{\log(m/2h)} {2 y^{h-1} \over \Gamma(h) m}   \, dy + n^{-\gamma} \int_{\log(m/2h)}^\infty {y^{h-1} e^{-y} \over \Gamma(h)}  \, dy \cr
&\le  {2 (\log(n))^h  \over m \, h! }  + { n^{-\gamma} (\log(n))^{h-1} 4h \over \Gamma(h) m}   \quad \hbox{\rm (for $n$ large enough, by lemma 1)} \cr
&= O \left(  n^{o(1)} / m  \right) 
= O \left(  m^{-1+o(1)}  \right)  . \cr}
$$
Thus, our probability (7) is not larger than
$O \left(  m^{ -2 + o(1)}  \right)$.
If there is only one unwanted parent selection and we avoid indices below $n^\gamma$,
and considering that the first parent selection at the root node is always
good, we see that at least half of the $2^h$ potential leaves are in fact
realized. Each of these leaves makes two independent parent selections. The probability
that all these leaves avoid parents $j$ with $S_j < (\sigma+\epsilon) \log n$ is at most
$(7/8)^{2^{h-1}} = o(n^{-2})$. If there is a connection, however, to such a parent of low depth,
then the root has shortest path length at most $h+1$ more than $(\sigma+\epsilon) \log n$.
Hence, if ${\cal E}_m$ is the event 
$\left[ S_m > (\sigma + \epsilon) \log n + h + 1 \right]$,
then
$$
\PROB \left\{ {\cal E}_m \cap ((B(r))^c \cap  \left[ \max_{1 \le j \le \lfloor n^{\gamma} \rfloor } S_j \le (\sigma + \epsilon) \log n \right] \right\}
= O \left( m^{-2 + o(1)} \right).
$$

Thus
$$
\eqalignno{
\PROB &\left\{ \max_{n^\gamma \le j \le n} S_j > (\sigma + \epsilon) \log n + h + 1 \right\} \cr
&= \PROB \left\{ \cup_{m\ge n^\gamma}^n  {\cal E}_m  \right\}  \cr
&\le 
\PROB \left\{ \max_{1 \le j \le \lfloor n^{\gamma} \rfloor } S_j > (\sigma + \epsilon) \log n \right\}
+
\PROB \{ B(r) \}
+
\sum_{m \ge n^\gamma}^n m^{-2 + o(1)} \cr
&= 
O(n^{-3\gamma }) + O(1/r') + n^{-\gamma +o(1)} \cr
&=
 O(1/\log n). \cr}
$$
This concludes the proof of the theorem. $\square$
\bigskip

\section
4. Bibliographic remarks and possible extensions.

The study of the {\sc urrt} goes back as far as
Na and Rapoport (1970) and Meir and Moon (1978).
Single nonuniform parent selections have been
considered as early as 1987 by Szyma\'nski.
Szyma\'nski (1987) showed that
if a parent is selected with probability proportional to its degree,
then with high probability
there is a node of  degree $\Omega (\sqrt{n} )$.
This is nothing but the preferential attachment model
of Barabasi and Albert (see Albert, Barabasi and Jeong, 1999,
or Albert and Barabasi, 1999), which for a single
parent is a special case of
the linear recursive trees  or {\sc port} (plane-oriented
recursive tree). 
For this model, the parameter $R_n$ was studied by Mahmoud (1992a),
and the height by Pittel (1994) and Biggins and Grey (1997), and
in a rather general setting by Broutin and Devroye (2006): the height
is in probability $(1.7956\ldots + o(1)) \log n$.
The profile (number of nodes at each depth level) was
studied by Hwang (2005, 2007) and Sulzbach (2008).

One can ask the questions studied in the present paper
for these more general models.

Various aspects of {\sc urrt}'s besides the depth and height
have been studied by many researchers. These include the degrees
of the nodes, the profile,
sizes of certain subtrees of certain nodes, the number of leaves,
and so forth. Surveys and references can be found in the book
by Mahmoud (1992b) or the paper by Devroye (1998).
Specific early papers include Timofeev (1984),
Gastwirth (1997),
Dondajewski and Szyma\'nski (1982),
Mahmoud (1991),
Mahmoud and Smythe (1991),
Smythe and Mahmoud (1994),
Szyma\'nski (1990), 
and the most recent contributions include
Fuchs, Hwang and Neininger (2006),
and
Drmota, Janson and Neininger (2008).
One may wonder how the profiles behave for 
uniform random $k$-dags.

\bigskip 
\section
4. References

\parskip=0pt 
\hyphenpenalty=-1000 \pretolerance=-1 \tolerance=1000 
\doublehyphendemerits=-100000 \finalhyphendemerits=-100000 
\frenchspacing 
\def\beginref{ 
\par\begingroup\nobreak\smallskip\parindent=0pt\kern1pt\nobreak 
\everypar{\strut}  } 
\def\endref{ 
\kern1pt\endgroup\smallbreak\noindent} 
\beginref

\endref

\beginref R.~Albert
and A.~Barabasi,
1999,
``Emergence of scaling in random networks,''
{\sl Science},
vol.~286,
pp.~509--512.

\endref

\beginref R.~Albert,
A.~Barabasi,
and H.~Jeong,
1999,
``Diameter of the World-Wide Web,''
{\sl Nature},
vol.~401,
p.~130.

\endref

\beginref J.~L.~Balcazar,
J.~Diaz,
and J.~Gabarro,
1995,
{\sl Structural Complexity I},
Springer-Verlag,
Berlin.

\endref

\beginref J.~D.~Biggins,
1976,
``The first and last-birth problems for a multitype age-dependent branching process,''
{\sl Advances in Applied Probability},
vol.~8,
pp.~446--459.

\endref

\beginref J.~D.~Biggins,
1977,
``Chernoff's theorem in the branching random walk,''
{\sl Journal of Applied Probability},
vol.~14,
pp.~630--636.

\endref

\beginref J.~D.~Biggins
and D.~R.~Grey,
1997,
``A note on the growth of random trees,''
{\sl Statistics and Probability letters},
vol.~32,
pp.~339--342.

\endref

\beginref N.~Broutin
and L.~Devroye,
2006,
``Large deviations for the weighted height of an extended class of trees,''
{\sl Algorithmica},
vol.~46,
pp.~271--297.

\endref

\beginref N.~Broutin,
L.~Devroye,
and E.~McLeish,
2008,
``Weighted height of random trees,''
{\sl Acta Informatica},
vol.~45,
pp.~237--277.

\endref

\beginref B.~Codenotti,
P.~Gemmell,
and J.~Simon,
1995,
``Average circuit depth and average communication complexity,''
in: {\sl Third European Symposium on Algorithms},
pp.~102--112.
Springer-Verlag,
Berlin.

\endref

\beginref L.~Devroye,
1986,
``A note on the height of binary search trees,''
{\sl Journal of the ACM},
vol.~33,
pp.~489--498.

\endref

\beginref L.~Devroye,
1987,
``Branching processes in the analysis of the heights of trees,''
{\sl Acta Informatica},
vol.~24,
pp.~277--298.

\endref

\beginref L.~Devroye,
1988,
``Applications of the theory of records in the study of random trees,''
{\sl Acta Informatica},
vol.~26,
pp.~123--130.

\endref

\beginref L.~Devroye,
1998,
``Branching processes and their applications in the analysis of tree structures and tree algorithms,''
in: {\sl Probabilistic Methods for Algorithmic Discrete Mathematics},
edited by M.~Habib, C.~McDiarmid, J.~Ramirez-Alfonsin and B.~Reed,
vol.~16,
pp.~249--314.
Springer-Verlag,
Berlin.

\endref

\beginref L.~Devroye,
1999,
``Universal limit laws for depths in random trees,''
{\sl SIAM Journal on Computing},
vol.~28,
pp.~409--432.

\endref

\beginref J.~Diaz,
M.~J.~Serna,
P.~Spirakis,
J.~Toran,
and T.~Tsukiji,
1994,
``On the expected depth of Boolean circuits,''
Technical Report LSI-94-7-R, Universitat Politecnica de Catalunya, Dep.\ LSI.

\endref

\beginref E.~W.~Dijkstra,
1959,
``A note on two problems in connexion with graphs,''
{\sl Numerische Mathematik},
vol.~1,
pp.~269--271.

\endref

\beginref M.~Dondajewski
and J.~Szyma\'nski,
1982,
``On the distribution of vertex-degrees in a strata of a random recursive tree,''
{\sl Bulletin de l'Acad\'emie Polonaise des Sciences, S\'erie des Sciences Math\'ematiques},
vol.~30,
pp.~205--209.

\endref

\beginref M.~Drmota,
S.~Janson,
and R.~Neininger,
2008,
``A functional limit theorem for the profile of search trees,''
{\sl Annals of Applied Probability},
vol.~18,
pp.~288--333.

\endref

\beginref M.~Fuchs,
H.-K.~Hwang,
and R.~Neininger,
2006,
``Profiles of random trees: Limit theorems for random recursive trees and binary search trees,''
{\sl Algorithmica},
vol.~46,
pp.~367--407.

\endref

\beginref D.~K.~Fuk
and S.~V.~Nagaev,
1971,
``Probability inequalities for sums of independent random variables,''
{\sl Theory of Probability and its Applications},
vol.~16,
pp.~643--660.

\endref

\beginref J.~L.~Gastwirth,
1977,
``A probability model of a pyramid scheme,''
{\sl The American Statistician},
vol.~31,
pp.~79--82.

\endref

\beginref N.~Glick,
1978,
``Breaking records and breaking boards,''
{\sl American Mathematical Monthly},
vol.~85,
pp.~2--26.

\endref

\beginref H.-K.~Hwang,
2005,
``Profiles of random trees: plane-oriented recursive trees (Extended Abstract),''
in: {\sl International Conference on Analysis of Algorithms,DMTCS Proceedings AD},
pp.~193--200.

\endref

\beginref H.-K.~Hwang,
2007,
``Profiles of random trees: Plane-oriented recursive trees,''
{\sl Random Structures and Algorithms},
vol.~30,
pp.~380--413.

\endref

\beginref H.~M.~Mahmoud,
1991,
``Limiting distributions for path lengths in recursive trees,''
{\sl Probability in the Engineering and Informational Sciencies},
vol.~5,
pp.~53--59.

\endref

\beginref H.~Mahmoud,
1992a,
``Distances in plane-oriented recursive trees,''
{\sl Journal of Computers and Applications in Mathematics},
vol.~41,
pp.~237--245.

\endref

\beginref H.~Mahmoud,
1992b,
{\sl Evolution of Random Search Trees},
Wiley,
New York.

\endref

\beginref H.~M.~Mahmoud
and R.~T.~Smythe,
1991,
``On the distribution of leaves in rooted subtrees of recursive trees,''
{\sl Annals of Applied Probability},
vol.~1,
pp.~406--418.

\endref

\beginref H.~Mahmoud
and B.~Pittel,
1984,
``On the most probable shape of a search tree grown from a random permutation,''
{\sl SIAM Journal on Algebraic and Discrete Methods},
vol.~5,
pp.~69--81.

\endref

\beginref H.~S.~Na
and A.~Rapoport,
1970,
``Distribution of nodes of a tree by degree,''
{\sl Mathematical Biosciences},
vol.~6,
pp.~313--329.

\endref

\beginref V.~V.~Petrov,
1995,
{\sl Limit Theorems of Probability Theory: Sequences of Independent Random Variables},
Clarendon Press,
Oxford.

\endref

\beginref B.~Pittel,
1984,
``On growing random binary trees,''
{\sl Journal of Mathematical Analysis and Applications},
vol.~103,
pp.~461--480.

\endref

\beginref B.~Pittel,
1994,
``Note on the heights of random recursive trees and random m-ary search trees,''
{\sl Random Structures and Algorithms},
vol.~5,
pp.~337--347.

\endref

\beginref R.~C.~Prim,
1957,
``Shortest connection networks and some generalizations,''
{\sl BSTJ},
vol.~36,
pp.~1389--1401.

\endref

\beginref R.~Pyke,
1965,
``Spacings,''
{\sl Journal of the Royal Statistical Society Series B},
vol.~7,
pp.~395--445.

\endref

\beginref A.~R\'enyi,
1962,
``Theorie des elements saillant d'une suite d'observations,''
in: {\sl Colloquium on Combinatorial Methods in Probability Theory},
pp.~104--115.
Mathematisk Institut,
Aarhus Universitet, Denmark.

\endref

\beginref H.~P.~Rosenthal,
1970,
``On the subspaces of $L^p$ ($p > 2$) spanned by sequences of independent random variables,''
{\sl Israel Journal of Mathematics},
vol.~8,
pp.~273--303.

\endref

\beginref R.~T.~Smythe
and H.~M.~Mahmoud,
1994,
``A survey of recursive trees,''
{\sl Teorya Imovirnostyta Mat. Stat. (in Ukrainian)},
vol.~51,
pp.~1--29.

\endref

\beginref H.~Sulzbach,
2008,
``A functional limit law for the profile of plane-oriented recursive trees.,''
in: {\sl Fifth Colloquium on Mathematics and Computer Science, DMTCS Proceedings AI},
pp.~339--350.

\endref

\beginref J.~Szyma\'nski,
1987,
``On a nonuniform random recursive tree,''
{\sl Annals of Discrete Mathematics},
vol.~33,
pp.~297--306.

\endref

\beginref J.~Szyma\'nski,
1990,
{\sl On the maximum degree and height of a random recursive tree},
Wiley,
New York.

\endref

\beginref E.~A.~Timofeev,
1984,
``Random minimal trees,''
{\sl Theory of Probability and its Applications},
vol.~29,
pp.~134--141.

\endref

\beginref T.~Tsukiji
and F.~Xhafa,
1996,
``On the depth of randomly generated circuits,''
in: {\sl Proceedings of Fourth European Symposium on Algorithms}.

\endref

\beginref 
\endref

\beginref 
\endref \bye